\documentclass[letterpaper]{article}
\usepackage{amsmath, amsthm, amssymb, diagrams}

\theoremstyle{plain}
\newtheorem{thm}{Theorem}
\newtheorem*{thm*}{Theorem}
\newtheorem{lemma}[thm]{Lemma}
\newtheorem{prop}[thm]{Proposition}
\newtheorem*{prop*}{Proposition}

\newtheorem*{question*}{Question}

\theoremstyle{definition}
\newtheorem*{definition*}{Definition}
\newtheorem*{example*}{Example}
\newtheorem{remark}[thm]{Remark}
\newtheorem*{remark*}{Remark}

\renewcommand{\hom}{\operatorname{Hom}_R}
\newcommand{\ext}{\operatorname{Ext}_R}
\newcommand{\pd}[1]{\operatorname{pd}\!\left(#1\right)}
\newcommand{\pdr}[2]{\operatorname{pd}_#1\!\left(#2\right)}
\newcommand{\codim}[1]{\operatorname{ht}\!\left(#1\right)}
\newcommand{\grade}[1]{\operatorname{grade}(#1)}
\newcommand{\depth}[1]{\operatorname{depth}(#1)}
\newcommand{\ass}[1]{\operatorname{Ass}\!\left(#1\right)}

\newcommand{\cm}{Cohen-Macaulay}
\newcommand{\unm}{^\mathrm{unm}}
\newcommand{\m}{\mathfrak{m}}
\newcommand{\fgh}{\ensuremath{(f,g,h)}}
\newcommand{\pdfgh}{\pd{R/\fgh}}
\def\sfrac#1#2{\hbox{$^{#1} \!/\! _{#2}$}}

\title{\bf On the projective dimension and the unmixed part of three cubics}

\author{Bahman Engheta
		\\
		\\
		\small\it
		Department of Mathematics, University of California, Riverside, CA 92521, USA
		\\
		\small\it
		E-mail address: \rm engheta@math.ucr.edu}
\date{}

\begin{document}

\maketitle 

\begin{abstract}
Let $R$ be a polynomial ring over a field in an unspecified number of variables. 
We prove that if $J \subset R$ is an ideal generated by three cubic forms, and the unmixed part of $J$ contains a quadric, then the projective dimension of $R/J$ is at most 4. 
To this end, we show that if $K \subset R$ is a three-generated ideal of height two and $L \subset R$ an ideal linked to the unmixed part of $K$, then the projective dimension of $R/K$ is bounded above by the projective dimension of $R/L$ plus one. 
\end{abstract}


\section{Introduction}
\label{section:intro}

Throughout this paper, unless otherwise stated, let $R$ denote \emph{any} polynomial ring over an arbitrary field $k$, say, $R=k[X_1, \ldots, X_n]$ where $n$ is not specified. This work was motivated by the following question posed by Michael E. Stillman. 

\begin{question*}[Stillman]
Is there a bound, independent of $n$, on the projective dimension of ideals in $R=k[X_1, \ldots, X_n]$ which are generated by $N$ homogeneous polynomials of given degrees $d_1,\ldots,d_N$?
\end{question*}

Equivalently, when $R$ denotes a polynomial ring over a field in an \emph{unknown number of variables}, does the module $F_1 \cong \bigoplus_{j=1}^N R(-d_j)$ in an arbitrary minimal graded free resolution of the form $~\cdots \to F_2 \to F_1 \to R~$ determine a bound on the length of this resolution?

This question concerns the existence of a uniform bound on the projective dimension of $R/J$ where neither the ring $R$ nor the ideal $J \subset R$ are fixed, but merely the number of generators of $J$ and the degrees of those generators. 
In other words, it asks whether the quantity 
\begin{multline*}
	\sup_n \big\{ \, \pdr{R}{R/J} ~~ | ~~ J \subset R=k[X_1,\ldots,X_n] \mbox{ is an ideal} \\
	\mbox{generated by } N \mbox{ forms of degrees } d_1, \ldots, d_N \, \big\}
\end{multline*}
is finite, where $\pdr{R}{R/J}$ denotes the projective dimension of $R/J$ over $R$. 
Henceforth we shall omit the subscript $R$ and write $\pd{R/J}$ for short.

\begin{remark*}[\it Projective dimension of three-generated ideals] 
At this juncture it is worth recalling the construction of Burch~\cite{burch,kohn}, whereby for any arbitrarily large integer $s$ one can construct a three-generated ideal $J_s \subset R=k[X_1,\ldots,X_{2s-4}]$ with $\pdr{R}{R/J_s}=s$. We note, however, that this fact does not furnish an answer to Stillman's question, for Burch's construction imposes lower bounds on the degrees of the generators of $J_s$ which increase with $s$.
\end{remark*}

Clearly, the above question has an affirmative answer when $N \leq 2$ or when $d_1=\cdots=d_N=1$. Already for $N=3$ and $d_1=d_2=d_3=2$, however, nontrivial arguments are needed to prove the existence of a bound. It was verified by Eisenbud and Huneke that indeed, if $J$ is generated by three quadric forms, then $\pd{R/J} \leq 4$. Considerably more effort is required to prove the existence of a bound on the projective dimension of three cubic forms; the author has shown in \cite{engheta-thesis} that if $J$ is generated by three cubic forms, then $\pd{R/J} \leq 36$. In this paper we present the following result which serves as a stepping stone in that direction. 

\begin{thm*}
Let $R$ be a polynomial ring over a field and let $J \subset R$ be an ideal generated by three cubic forms. 
If the unmixed part of $J$ contains a quadric, then $\pd{R/J} \leq 4$. 
\end{thm*}

Recall that the unmixed part of an ideal in a Noetherian ring is the intersection of its primary components of minimal height. Our approach, which concentrates on the unmixed part of the ideal $J$ and involves linkage theory, was motivated by the works of Huneke and Ulrich~\cite{hu2}, and Fan~\cite{fan}. 
Building on their result, we also exhibit the following bound which does not require any assumptions on the degrees of the generators of the ideal. 

\begin{prop*}
Let $R$ be a polynomial ring over a field and let $J \subset R$ be a three-generated ideal of height at least two. 
If the unmixed part of $J$ contains a linear form, then $\pd{R/J} \leq 3$. 
\end{prop*}

This paper is organized as follows. 
In the remainder of this section we reduce the question concerning a bound on the projective dimension of an ideal generated by three cubic forms to the case where the ideal has height two and the ground field $k$ is algebraically closed. 
In Section~\ref{section:unmix-link} we recall the notion of algebraic linkage and prove the following. 

\begin{thm*}
Let $R$ be a regular local ring, let $J \subset R$ be a three-generated ideal of height two, and let $I$ denote the unmixed part of $J$. 
If $I'$ is an ideal linked to $I$, then $\pd{R/J} \leq \pd{R/I'}+1$. 
\end{thm*}

In Section~\ref{section:low-mult} we study height two unmixed ideals in $R$ of multiplicity $\leq 3$. 
We classify all such ideals of multiplicity 2 and show that any three-generated ideal in $R$ of height two and multiplicity 2 has projective dimension $\leq 4$. 
In Section~\ref{section:bounds} we show that if $J$ is generated by three cubic forms with multiplicity 3, then $\pd{R/J} \leq 16$. 
We conclude that section by proving the aforementioned result and the following theorem which is similar in nature. 

\begin{thm*}
Let $R$ be a polynomial ring over a field and let $J \subset R$ be an ideal generated by three cubic forms with multiplicity $\leq 5$. 
If an ideal which is linked to the unmixed part of $J$ contains a quadric, then $\pd{R/J} \leq 4$. 
\end{thm*}

\subsection{Preliminaries}
\label{subsection:prelims}

{\it Notation.} 
We will denote by $\mathfrak m$ the homogeneous maximal ideal $(X_1,\ldots,X_n)$ of $R$. 
For an ideal $J \subset R$, $\codim{J}$ denotes the height of $J$ and $J\unm$ is the unmixed part (or the top dimensional component) of $J$, that is, the intersection of those primary components $Q$ of $J$ with $\codim{Q} = \codim{J}$. 
Note that in general $\sqrt{J} \subsetneq J\unm$. 
By $\lambda(R/J)$ we denote the length of $R/J$ and by $e(R/J)$ its multiplicity at $\mathfrak m$. 
One has $e(R/J) = e(R/J\unm)$, as easily follows from the associativity formula for multiplicities: 
\begin{equation}
\label{eqn:assoc}
	e(R/J) ~ = \!\!\!\!\! 
	\sum_{\substack{P \, \in \, \mathrm{Spec}(R) \\ \dim(R/P) \, = \, \dim(R/J) }} 
	\!\!\!\!\! 
	e(R/P) \, \lambda(R_P/J_P).
\end{equation}

Throughout our arguments, we will often employ the following well known lemma whenever we encounter short exact sequences. 

\begin{lemma}[Depth Lemma]
\label{lemma:depth}
Let $0\to A\to B\to C\to 0$ be a short exact sequence of finite $R$-modules. 
The following inequalities hold. 
\begin{align*}
	\pd{A} & \leq \max\{ \pd{B}, \, \pd{C}-1 \}, \\
	\pd{B} & \leq \max\{ \pd{A}, \, \pd{C} \}, \\
	\pd{C} & \leq \max\{ \pd{A}+1, \, \pd{B} \}.
\end{align*}
Furthermore, one has 
\begin{align*}
\label{eqn:more-depth}
	\pd{B}	\leq \pd{C}-1& ~~ \Longrightarrow ~~ \pd{A} = \pd{C}-1, \\
	\pd{C}	\leq \pd{A} ~~~~\, & ~~ \Longrightarrow ~~ \pd{B} = \pd{A}, \\
	\pd{A}+1	\leq \pd{B} ~~~~\, & ~~ \Longrightarrow ~~ \pd{C} = \pd{B}. 
\end{align*}
\end{lemma}

\begin{remark}[\it Reduction to the height two case]
\label{rmk:codim-2} 
The question concerning a bound on the projective dimension of three cubic forms $f,g,h$ can be reduced to the case where the three cubics generate an ideal of height two. 
Indeed, $\codim{f,g,h} \leq 3$ by Krull's Height Theorem. 
If $\codim{f,g,h}=1$, then $\fgh$ is contained in a prime ideal of height one. 
As $R$ is a unique factorization domain, this prime ideal is principal. 
So the cubics $f,g,h$ share a common divisor and $\fgh$ is isomorphic to an ideal generated either by three quadrics, in which case $\pdfgh \leq 4$, or by three linear forms, in which case $\pdfgh \leq 3$. 
And if $\codim{f,g,h}=3$, then $f,g,h$ form a regular sequence and $\pdfgh=3$. 
So one only needs to consider the case $\codim{f,g,h}=2$. 
\end{remark}

\begin{remark}[\it Reduction to an algebraically closed field]
\label{rmk:alg-closed} 
We may assume without loss of generality that the field $k$ is algebraically closed. 
For if $\bar k$ denotes the algebraic closure of $k$, then the map $R=k[X_1,\ldots,X_n] \longrightarrow S=\bar k[X_1,\ldots,X_n]$ is a flat ring homomorphism and thus $\pdr{R}{R/J} = \pdr{S}{S/JS}$. 
In particular, we may assume that $k$ is infinite which gives us the leverage to apply ``prime avoidance.'' 
\end{remark}

\section{The unmixed part and linkage}
\label{section:unmix-link}

Recall the notion of \emph{algebraic linkage} as introduced by Peskine and Szpiro~\cite{ps}. 

\begin{definition*}
Two proper ideals $A$ and $B$ of height $g$ in a \cm{} ring $S$ are said to be \emph{(directly) linked} if there exists a maximal regular sequence $\mathbf z = z_1,\ldots,z_g$ in $A \cap B$ such that $A=(\mathbf z) : B$ and $B=(\mathbf z) : A$. 
\end{definition*}

The ideals $A$ and $B$ in the above definition are necessarily unmixed of height $g$ and their multiplicities are complementary to each other in the sense that $e( S/(\mathbf z) ) = e(S/A) + e(S/B)$. 
If the underlying ring is Gorenstein, then the unmixedness property of an ideal is also sufficient for that ideal to be linked to another ideal. 
More precisely, the following fundamental result of linkage theory asserts that in a Gorenstein ring one can always produce a link to an unmixed ideal and that the \cm{} property is preserved by this process. 

\begin{prop}[Peskine-Szpiro \cite{ps}, {\cite[Proposition 2.5]{hu1}}]
\label{prop:linkage}
Let $A$ be an unmixed ideal of height $g$ in a Gorenstein ring $S$. 
Let $\mathbf z = z_1,\ldots,z_g$ be a maximal regular sequence inside $A$ with $(\mathbf z) \neq A$ and set $B = (\mathbf z):A$. 
Then 
\begin{itemize}
	\item[(a)] $A = (\mathbf z):B$, that is, $A$ and $B$ are algebraically linked. 
	\item[(b)] $S/A$ is \cm{} if and only if $S/B$ is \cm{}. 
\end{itemize}
\end{prop}

Lastly we recall that any two links of an ideal in a Gorenstein ring have the same (finite or infinite) projective dimension. 

In the context of \emph{residual intersections}, a generalization of linkage theory, Huneke and Ulrich~\cite{hu2} deduced the following which was later also derived by Fan~\cite{fan} using homological algebra and which we shall extend in Theorem~\ref{thm:extend}. 

\begin{thm}[Huneke-Ulrich {\cite[Page 16]{hu2}}, Fan {\cite[Corollary 1.2]{fan}}]
\label{thm:fan}
Let $R$ be a regular local ring and let $J \subset R$ be a three-generated ideal of height two. Let $I$ denote the unmixed part of $J$. If $R/I$ is \cm{}, then $\pd{R/J} \leq 3$.
\end{thm}

By virtue of Theorem~\ref{thm:fan}, in our effort to bound the projective dimension of a three-generated ideal, we may restrict our attention to those ideals with a non-degenerate unmixed part: 

\begin{prop}
\label{prop:degen}
Let $R$ be a polynomial ring over a field and let $J \subset R$ be a three-generated ideal of height two. 
If the unmixed part of $J$ contains a linear form, then $\pd{R/J} \leq 3$.
\end{prop}
\begin{proof}
Let $I$ denote the unmixed part of $J$ and let $l \in I$ be a linear form.
Then $I/(l)$ is a height one, unmixed ideal in the unique factorization domain $R/(l)$, and therefore principal.
Lifting a generator of $I/(l)$ back to $R$, along with $l$, gives a generating set for the ideal $I$.
As $\codim{I}=2$, $R/I$ is a complete intersection and $\pd{R/J} \leq 3$ by Theorem~\ref{thm:fan}.
\end{proof}

The hypothesis in Theorem~\ref{thm:fan} that $R/I$ be \cm{} could also be stated as $\pd{R/I}=2$. 
In the context of Theorem~\ref{thm:fan}, Huneke asked whether $\pd{R/I}=t$ would imply $\pd{R/J} \leq t+1$. 
This was answered in the negative in~\cite[Example 1.4]{fan}. 
We seek to extend Theorem~\ref{thm:fan} by interpreting the \cm{} assumption on $R/I$ from the point of view of linkage. 
By part~(b) of Proposition~\ref{prop:linkage}, if $R/I$ is \cm{}, then so is $R/I'$ for any ideal $I'$ which is linked to $I$. 
In particular, $\pd{R/I'}=2$ and the conclusion of Theorem~\ref{thm:fan} could be stated as $\pd{R/J} \leq \pd{R/I'} + 1$. 
We show that this is the case in general, regardless of the value of $\pd{R/I'}$. 

\begin{thm}
\label{thm:extend}
Let $R$ be a regular local ring and let $J \subset R$ be a three-generated ideal of height two. 
Denote by $I$ the unmixed part of $J$ and let $a,b \in I$ be a regular sequence. 
Then 
\[
	\pd{R/J} \leq \pd{\sfrac{R}{(a,b):I}} + 1
\]
and equality holds if $\pd{R/J} \geq 4$.
\end{thm}
\begin{proof}
Let $M$ and $N$ denote second syzygy modules of $R/J$ and $R/I$, respectively. 
That is, there are short exact sequences
\begin{align}
	& 0 \to M \to F \to J \to 0, 
	\label{MFJ}
	\\
	& 0 \to N \to G \to I \to 0, 
	\label{NGI}
\end{align}
with free modules $F \cong R^3$ and $G$. 
We first exhibit a bound for the projective dimension of $M^* = \hom(M,R)$ in terms of $\pd{\sfrac{R}{(a,b):I}}$. 
Then, by drawing on a result of Bruns~\cite{bruns} on (oriented) second syzygy modules, we establish a bound for the projective dimension of $R/J$. 

In the proof of~\cite[Theorem 1.1]{fan}, Fan establishes a short exact sequence
\begin{equation}
\label{eqn:dagger}
	0 \to G^* \to F^* \oplus N^* \to M^* \to 0.
\end{equation}
As $F^*$ and $G^*$ are free modules of projective dimension 0, this short exact sequence implies $\pd{M^*} \leq \max\{ 1, \, \pd{N^*} \}$. 
Next we bound the projective dimension of $N^*$. 

As $\grade{I} = 2$, we have $\ext^1(R/I, R) = 0$ and the long exact sequence on $\ext^\bullet(-,R)$ induced by the short exact sequence $0\to I\to R\to R/I\to 0$ yields $I^* \cong R$. 
The same long exact sequence also yields $\ext^1(I,R) \cong \ext^2(R/I,R)$. 
So dualizing \eqref{NGI} gives rise to the exact sequences
\begin{diagram}[size=1.1em]
	&0\rTo	&R	&\rTo&G^*&		&\rTo&&N^*&\rTo&\ext^2(R/I,R)&\rTo&0, \\
	&		&	&	&	&\rdTo	&	&\ruTo \\
	&		&	&	&	&		&K\\
	&		&	&	&	&\ruTo	&	&\rdTo \\
	&		&	&	&0	&		&	&&0
\end{diagram}
where $K$ is the image of $G^* \to N^*$. 
As $0 \to R \to G^* \to K \to 0$ is a free resolution of $K$, we have $\pd{K} \leq 1$. 
Thus, $\pd{N^*} \leq \max\!\left\{ 1, \, \pd{\ext^2(R/I,R)} \right\}$. 

In order to bound the projective dimension of $\ext^2(R/I,R)$, recall that for any regular sequence $a,b \in I$ one has 
\[
	\ext^2(R/I,R) \cong \hom(R/I,R/(a,b)) \cong \sfrac{(a,b):I}{(a,b)}, 
\]
cf.~\cite[Lemma 1.2.4]{bh}. 
And the short exact sequence 
\begin{equation}
\label{eqn:link-ses}
	0\longrightarrow \frac{(a,b):I}{(a,b)}\longrightarrow \frac{R}{(a,b)}\longrightarrow \frac{R}{(a,b):I}\longrightarrow 0
\end{equation}
delivers $\pd{\sfrac{(a,b):I}{(a,b)}} \leq \max\!\left\{ 2, \, \pd{\sfrac{R}{(a,b):I}} - 1 \right\}$. 
Combining the inequalities obtained so far, we arrive at $\pd{M^*} \leq \max\!\left\{ 2, \, \pd{\sfrac{R}{(a,b):I}} - 1 \right\}$. 

As $J$ is a three-generated ideal, the free module $F$ in the short exact sequence \eqref{MFJ} has rank three, whence $M$ is a second syzygy module of rank two. 
A result of Bruns~\cite[Corollary 2.6]{bruns} now states that $M \cong M^*$. 
So $\pd{R/J} = \pd{M}+2 = \pd{M^*}+2$ and we arrive at 
\begin{equation}
\label{eqn:max}
	\pd{R/J} \leq \max\!\left\{ 4, \, \pd{\sfrac{R}{(a,b):I}} + 1 \right\}.
\end{equation}

We have $\pd{\sfrac{R}{(a,b):I}} \geq \grade{(a,b):I} = 2$. 
If $\pd{\sfrac{R}{(a,b):I}} = 2$, then $\sfrac{R}{(a,b):I}$ is \cm{} and therefore so is $R/I$ by Proposition~\ref{prop:linkage}. 
In this case Theorem~\ref{thm:fan} asserts that $\pd{R/J} \leq 3 = \pd{\sfrac{R}{(a,b):I}} + 1$, as claimed. 
And if $\pd{\sfrac{R}{(a,b):I}} > 2$, then the desired inequality $\pd{R/J} \leq \pd{\sfrac{R}{(a,b):I}} + 1$ follows directly from \eqref{eqn:max}. 

It remains to show that $\pd{R/J} = \pd{\sfrac{R}{(a,b):I}} + 1$ whenever $\pd{R/J} \geq 4$. 
Set $j = \pd{R/J}$ and assume $j \geq 4$. 
So we have $\pd{M^*} = \pd{M} = j-2 \geq 2$. 
It follows from the short exact sequence \eqref{eqn:dagger} that $\pd{N^*} = \pd{M^*}$. 
As $\pd{K} \leq 1$ and $\pd{N^*} \geq 2$, the short exact sequence 
\[
	0 \to K \to N^* \to \sfrac{(a,b):I}{(a,b)} \to 0
\]
implies that $\pd{\sfrac{(a,b):I}{(a,b)}} = \pd{N^*} = j-2$. 
Finally, by the short exact sequence \eqref{eqn:link-ses}, $\pd{\sfrac{R}{(a,b):I}} \leq \max\{ j-1, \, 2 \}$. 
As we are assuming $j = \pd{R/J} \geq 4$, this maximum equals $j-1$. 
That is, $\pd{R/J} \geq \pd{\sfrac{R}{(a,b):I}} + 1$. 
This finishes the proof, as we have already established the reverse inequality above. 
\end{proof}

Theorem~\ref{thm:extend} proves to be useful even in instances where one cannot determine the unmixed part $I$ explicitly, but where one can choose elements $a,b \in I$ of sufficiently low degree in order to bound the multiplicity of $I$ and of its link $(a,b):I$, and use this information on multiplicity to give an upper bound for $\pd{\sfrac{R}{(a,b):I}}$ and consequently for $\pd{R/J}$. 
See proof of Theorem~\ref{thm:unmix-quad} for such an application.

\section{Unmixed ideals of low multiplicity}
\label{section:low-mult}

The aim of this section is to establish properties of height two unmixed ideals of multiplicity 2 and 3 which will be used in the proofs of the subsequent results. 
We begin with two elementary lemmata. 

\begin{lemma}
\label{lemma:unmix-eq}
Let $J \subset R$ be an unmixed ideal. 
If $I \subset R$ is an ideal containing $J$, with the same height and multiplicity as $J$, then $J=I$.
\end{lemma}
\begin{proof}
As $J\subseteq I$ and $\codim{J} = \codim{I}$, the prime ideals contributing to the multiplicity of $R/I$ also contribute to the multiplicity of $R/J$, that is, 
\begin{equation}
\label{eqn:contrib}
	\{ P \in \ass{R/I} ~|~ \codim{P}=\codim{I} \} \subseteq \{ P \in \ass{R/J} ~|~ \codim{P}=\codim{J} \}.
\end{equation}
Also, the inclusion $J \subseteq I$ implies $\lambda(R_P/J_P) \geq \lambda(R_P/I_P)$ for $P \in \ass{R/I}$.
So we have the following two inequalities: 
\begin{equation*} \begin{split}
	e(R/J) ~ & = \sum_{\substack{P \, \in \, \ass{R/J} \\ \codim{P} =\codim{J} }} 
	e(R/P) \; \lambda(R_P/J_P) \\
	& \geq \sum_{\substack{P \, \in \, \ass{R/I} \\ \codim{P}=\codim{I} }} 
	e(R/P) \; \lambda(R_P/J_P) \\
	& \geq \sum_{\substack{P \, \in \, \ass{R/I} \\ \codim{P}=\codim{I} }} 
	e(R/P) \; \lambda(R_P/I_P) ~ = ~ e(R/I).
\end{split} \end{equation*}
As $e(R/J)=e(R/I)$ by our hypothesis, this entails that equality holds in \eqref{eqn:contrib} and $\lambda(R_P/J_P) = \lambda(R_P/I_P)$ for all $P\in\ass{R/J}$ with $\codim{P} = \codim{J} $. 
Since $J$ is unmixed, these constitute all the prime ideals associated to $J$.
So $J_P=I_P$ for all $P\in\ass{R/J}$ and therefore $J=I$.
\end{proof}

\begin{lemma}
\label{lemma:unmixed-hull}
If $J \subseteq I$ are two ideals with $\codim{J}=\codim{I}$, then $J\unm \subseteq I\unm$.
\end{lemma}
(The assumption on height is necessary; 
e.g.\ consider the ideals $J=(X) \cap (Y,Z)$ with $J\unm = (X)$ and $I=I\unm=(Y,Z)$.) 
\begin{proof}
The statement of the lemma is a consequence of the following more general fact. 
Every primary component of $I$ of minimal height contains the corresponding primary component of $J$. 
More precisely, if $P \in \ass{R/I}$ with $\codim{P}=\codim{I}$, then $P \in \ass{R/J}$ with $\codim{P}=\codim{J}$ by our hypotheses. 
So if $I\unm = K_1 \cap \cdots \cap K_m$ with $\sqrt{K_i} = P_i$, then $J\unm = L_1 \cap \cdots \cap L_m \cap L_{m+1} \cap \cdots \cap L_{n}$ with $\sqrt{L_i} = P_i$ for $i=1 \ldots m$. 
Now, showing $L_i \subseteq K_i$ for $i=1 \ldots m$ is sufficient to prove the inclusion $J\unm \subseteq I\unm$. 

Since $P_i$ is a minimal prime of $J$ and $I$, localizing $J \subseteq I$ at $P_i$ yields $(L_i)_{P_i} \subseteq (K_i)_{P_i}$. 
As $L_i$ and $K_i$ are primary to $P_i$, this entails $L_i \subseteq K_i$.
\end{proof}

The following lemma describes a class of unmixed (in fact, primary) ideals which we will often encounter. 
Note that locally at their associated prime $(x,y)$, and after a linear change of coordinates, these ideals are simply of the form $(x,y^e)_{(x,y)}$, where $e$ is a positive integer. 

\begin{lemma}
\label{lemma:unmix}
Let $R$ be a polynomial ring over a field and let $I = (x,y)^e + (ax+by)$ with independent linear forms $x,y$, a positive integer $e$, and forms $a,b \in R$ (of equal degree) such that $(a,b) \not\subset (x,y)$. 
Then $\pd{R/I} \leq 3$. 
Further, $I$ is unmixed if and only if $\codim{x,y,a,b}>3$, in which case $I$ is primary to $(x,y)$ and $R/I$ has multiplicity $e$. 
\end{lemma}
\begin{proof}
Consider the sequence of free $R$-modules 
\begin{equation}
\label{eqn:complex}
	0 \to R^{e-1}
	\xrightarrow{\varphi_3} R^{2e}
	\xrightarrow{\varphi_2} R^{e+2}
	\xrightarrow{\varphi_1} R
	\to 0, 
\end{equation}
with the maps 
\[
	\varphi_1 = \big( 
	\underbrace{ \begin{array}{lccccr}
		ax+by ~ & x^e ~ & x^{e-1}y ~ & \cdots ~ & xy^{e-1} ~ & y^e
	\end{array} }_{e+2 \mathrm{~columns}}
	\big),
\]

\[
\varphi_2 = 
\left( \begin{matrix}\\ \\ \\ \\ \\ \\ \\ \end{matrix} \right.
\underbrace{ \begin{array}{ccccc}
	x^{e-1} & x^{e-2}y & \cdots & xy^{e-2} & y^{e-1} 
	\\ \hline
	-a&  &		&  &
	\\
	-b&-a&		&  &
	\\
	  &-b&\ddots&  &
	\\
	  &  &\ddots&-a&
	\\
	  &  &		&-b&-a
	\\
	  &  &		&  &-b
\end{array} }_{e \mathrm{~columns}}
\!\!\!\!\!
\begin{array}{l|l}\\ \\ \\ \\ \\ \\ \\ \end{array}
\underbrace{ \begin{array}{rrrrr}
	&&&&
	\\ \hline
	 y&  &		&  &
	\\
	-x& y&		&  &
	\\
	  &-x&\ddots&  &
	\\
	  &  &\ddots& y&
	\\
	  &  &		&-x& y
	\\
	  &  &		&  &-x
\end{array} }_{e \mathrm{~columns}}
\left. \begin{matrix}\\ \\ \\ \\ \\ \\ \\ \end{matrix} \right)
\left. \begin{matrix}\\ \\ \\ \\ \\ \\ \\ \end{matrix} \right\} 
e+2 \mathrm{~rows,}
\]
and 
\[
\varphi_3 = 
\left( \begin{matrix}\\ \\ \\ \\ \\ \\ \\ \\ \\ \end{matrix} \right.
\underbrace{ \begin{array}{rrr}
	 y&      &
	\\
	-x&\ddots&
	\\
	  &\ddots& y
	\\
	  &      &-x
	\\ \hline
	 a&      &
	\\
	 b&\ddots&
	\\
	  &\ddots& a
	\\
	  &      & b
\end{array} }_{e-1 \mathrm{~columns}}
\left. \begin{matrix}\\ \\ \\ \\ \\ \\ \\ \\ \\ \end{matrix} \right)
\begin{matrix}
	\left. \begin{matrix}\\ \\ \\ \\ \end{matrix} \right\} 
	e \mathrm{~rows}
	\\
	\left. \begin{matrix}\\ \\ \\ \\ \end{matrix} \right\} 
	e \mathrm{~rows}
\end{matrix}
\]
(The missing entries of $\varphi_2$ and $\varphi_3$ are understood to be zero.) 
The compositions $\varphi_1 \varphi_2$ and $\varphi_2 \varphi_3$ are both zero, that is, \eqref{eqn:complex} is a complex of free $R$-modules. 
In fact, this complex is obtained as a \emph{generic} free resolution of $R/I$, meaning that the elements $a,b$ were chosen in such a way that $x,y,a,b$ form a regular sequence. 
As it turns out, \eqref{eqn:complex} is exact regardless of the choice of $a,b$ as long as $(a,b) \not\subset (x,y)$, which we shall show next. 

Let $I_{r_i}(\varphi_i)$ denote the ideal generated by the $r_i \times r_i$ minors of $\varphi_i$, where $r_i$ is the expected rank of $\varphi_i$. 
We have $r_1=1$, $r_2=e+1$, and $r_3=e-1$. 
One can verify by calculation that 
\[
	I_1(\varphi_1) = I, 
	\qquad
	I_{e+1}(\varphi_2) = I \, (x,y,a,b)^{e-1}, 
	\qquad
	I_{e-1}(\varphi_3) = (x,y,a,b)^{e-1}.
\]
(For the purposes of our proof it suffices to merely observe the inclusions $(x^{2e-1}, y^{2e-1}) \subset I_{e+1}(\varphi_2)$ and $(x^{e-1}, y^{e-1}, a^{e-1}, b^{e-1}) \subset I_{e-1}(\varphi_3)$, both of which are evident from the structure of the matrices.) 
Since by assumption $(x,y)$ is a prime ideal of height two and $(a,b) \not\subset (x,y)$, we have $\codim{x,y,a,b} \geq 3$ and consequently $\codim{I_{e-1}(\varphi_3)} \geq 3$. 
And clearly the ideals $I_{e+1}(\varphi_2)$ and $I_1(\varphi_1)$ both have height two. 
Therefore, by the Buchsbaum-Eisenbud acyclicity criterion, the complex \eqref{eqn:complex} is in fact a free resolution of $R/I$ and $\pd{R/I} \leq 3$. 
(More precisely, $\pd{R/I}=3$ if $a,b \in \m$ and $\pd{R/I}=2$ otherwise.) 

Let $P\in\ass{R/I}$ so that $\depth{R_P/I_P} = 0$.
By the Auslander-Buchsbaum formula, $\pd{R_P/I_P} = \codim{P}$ and therefore $\codim{P} \leq 3$. 
Now if $\codim{x,y,a,b}>3$, then $I_{e-1}(\varphi_3) \not\subset P$, while $I_e(\varphi_3)=0$. 
So the homomorphism $\varphi_3$ splits locally at $P$ and $\pd{R_P/I_P} = \codim{P} = 2$. 
That is, $I$ is unmixed of height two.

Conversely, if $\codim{x,y,a,b} = 3$, then $a$ and $b$ share a common divisor $c \in \m$ modulo $(x,y)$. 
Write $a \equiv ca'$ and $b \equiv cb'$ modulo $(x,y)$ and note that the hypothesis $(a,b) \not\subset (x,y)$ forces $c \notin (x,y)$ and $(a',b') \not\subset (x,y)$. 
In particular, $\codim{x,y,c} = 3$. 
We have $ax+by \equiv c\,(a'x+b'y)$ modulo $(x,y)^2$. 
Multiplying $ax+by$ with $(x,y)^{e-2}$ and reducing modulo $(x,y)^e$, we obtain $c\,(a'x+b'y)\,(x,y)^{e-2} \subset I$. 
So $(x,y,c) \subseteq I:(a'x+b'y)(x,y)^{e-2}$ which implies that $(x,y,c)$ is contained in an associated prime of $R/I$.
As $I$ has height two, it cannot be unmixed.

Finally, as $\sqrt I = (x,y)$, $I$ has only $(x,y)$ as its minimal prime and it is unmixed if and only if it is primary to $(x,y)$. 
And since $(x,y)^e \subset I$ and $(a,b) \not\subset (x,y)$, locally at $(x,y)$ the Hilbert function of $(R/I)_{(x,y)}$ is given by $(\underbrace{1,\ldots,1}_{e \mathrm{~times}})$ and thus $e(R/I)=e$.
\end{proof}

In the following proposition we classify all height two unmixed ideals of multiplicity 2. 
The significance of the proposition lies in part (iv). 

\begin{prop}
\label{prop:class}
Let $R$ be a polynomial ring over a field and let $I \subset R$ be a height two unmixed ideal of multiplicity 2. 
Then $\pd{R/I} \leq 3$ and $I$ is one of the following ideals. 
\begin{itemize}
	\item[(i)] 
	A prime ideal generated by a linear form and an irreducible quadric. 
	\item[(ii)] 
	$(x,y) \cap (x,v)=(x,yv)$ with independent linear forms $x,y,v$. 
	\item[(iii)] 
	$(x,y) \cap (u,v)=(xu,xv,yu,yv)$ with independent linear forms $x,y,u,v$. 
	\item[(iv)] 
	The $(x,y)$-primary ideal $(x,y)^2 + (ax+by)$ with independent linear forms $x,y$ and forms $a,b \in \m$ such that $x,y,a,b$ form a regular sequence. 
	\item[(iv$^\circ$)] 
	$(x,y^2)$ with independent linear forms $x,y$. 
\end{itemize}
\end{prop}
\begin{proof}
By the associativity formula for multiplicities, $I$ is either the intersection of two ideals each generated by two independent linear forms, in which case $I$ is of type (ii) or (iii), or $I$ is primary. 

If $I$ is not a prime ideal as in part (i), then it is is primary to a height two prime ideal $P$ of multiplicity one, say $P=(x,y)$ with independent linear forms $x,y$. As $\lambda(R_P/I_P)=2$, we have $P^2 \subsetneq I \subsetneq P$ and $I$ is generated by $P^2$ plus additional terms of the form $(a_ix+b_iy)$ with $(a_i,b_i) \not\subset (x,y)$. 
We claim that $I$ contains only one such term as a minimal generator, that is, $I=(x,y)^2+(ax+by)$ with $(a,b) \not\subset (x,y)$. 
To prove this, we choose one of the terms $a_ix+b_iy$ among the minimal generators of $I$, say $ax+by$, and first show that $\codim{x,y,a,b}>3$. 
(That is, either $a$ or $b$ is a unit or $x,y,a,b$ form a regular sequence.)
Indeed, as $(x,y)$ is prime of height two and $(a,b) \not\subset (x,y)$, the ideal $(x,y,a,b)$ has height at least 3. 
If $\codim{x,y,a,b} = 3$, then $a$ and $b$ have a common divisor $c \in \m$ modulo $(x,y)$. 
Writing $a \equiv ca'$ and $b \equiv cb'$ modulo $(x,y)$, we have $ax+by \equiv c(a'x+b'y)$ modulo $(x,y)^2$. 
As $(x,y)^2 \subset I$, the element $c(a'x+b'y)$ is a minimal generator of $I$ and $c$ is a zero-divisor on $R/I$. 
Since $I$ is primary to $(x,y)$, we must have $c \in (x,y)$. 
However, the condition $(a,b) \not\subset (x,y)$ implies that $c \notin (x,y)$ --- a contradiction. 
Thus, $\codim{x,y,a,b}>3$. 

Our claim now follows from Lemma~\ref{lemma:unmix} which establishes that the ideal $(x,y)^2+(ax+by)$ is unmixed, and Lemma~\ref{lemma:unmix-eq} which implies that $(x,y)^2+(ax+by)$ equals $I$.
If $a,b \in \m$, then $\codim{x,y,a,b} = 4$ and $I$ is of type (iv).
And if either $a$ or $b$ is a unit, then (after a linear change of coordinates) $I$ is of type (iv$^\circ$).

To finish the proof, we need to verify that the projective dimension of $R/I$ is at most 3. 
The ideals of type (i), (ii), and (iv$^\circ$) are complete intersections and in those cases $\pd{R/I} = 2$. 
As for part~(iii), applying Lemma~\ref{lemma:depth} to the short exact sequence 
\[
	0 \longrightarrow \frac{R}{I} \longrightarrow \frac{R}{(x,y)} \oplus \frac{R}{(u,v)} \longrightarrow \frac{R}{(x,y,u,v)} \longrightarrow 0
\]
yields $\pd{R/I}=3$. 
And in part~(iv) we have $\pd{R/I} = 3$ by Lemma~\ref{lemma:unmix}. 
\end{proof}

We apply Proposition~\ref{prop:class} in conjunction with Theorem~\ref{thm:extend} to point out the following fact which will later be used in the proof of Theorem~\ref{thm:unmix-quad}. 

Let $J \subset R$ be a three-generated ideal of height two. 
If $e(R/J) \leq 2$, then $\pd{R/J} \leq 4$. 
Indeed, let $I=J\unm$ denote the unmixed part of $J$. 
If $e(R/J)=1$, then $I$ is generated by linear forms and $\pd{R/J} \leq 3$ by Proposition~\ref{prop:degen}. 
So suppose $e(R/J)=2$ and $I$ does not contain a linear form. 
Then, by Proposition~\ref{prop:class}, either $I=(xu, xv, yu, yv)$ with independent linear forms $x, y, u, v$, or $I = (x, y)^2 + (ax+by)$ where $x, y, a, b \in \m$ form a regular sequence. 
In the former case we compute the link $(xu, yv) : I = (xu, yv, xy, uv)$ and see that $\pd{\sfrac{R}{(xu, yv, xy, uv)}} = 3$. 
Similarly, in the latter case we compute the link $(x^2, y^2) : I = (x, y)^2 + (ax-by)$ and note that $\pd{\sfrac{R}{(x, y)^2 + (ax-by)}} = 3$ by Lemma~\ref{lemma:unmix}. 
Thus, in both cases $\pd{R/J} \leq 4$ by Theorem~\ref{thm:extend}. 

\medskip

A classification similar to Proposition~\ref{prop:class} for height two unmixed ideals of multiplicity 3 remains elusive. 
This is due to the difficulty of determining all such primary ideals, as it was done in part (iv) of Proposition~\ref{prop:class}. 
Mimicking the proof of Proposition~\ref{prop:class} in the multiplicity 3 case leads to the following nontrivial example of such an ideal. 
(Trivial examples would be $(x,y)^2$ or $(x,y)^3+(cx+dy)$ where $x,y$ are independent linear forms and $x,y,c,d$ form a regular sequence.) 

\begin{example*}[\it A triple structure] 
Let $R=k[a,b,c,d,e,x,y]$ be a polynomial ring over a field $k$ and suppose $I \subset R$ is an ideal of multiplicity 3 and primary to $(x,y)$. 
In particular, $\lambda\!\left( (R/I)_{(x,y)} \right) = 3$ and $(x,y)^3 \subset I$. 
Set 
\begin{align*}
	q	= \ & (ac+dx)x + (bc+ey)y \\
		= \ & (ax+by)c + dx^2 + ey^2
\end{align*}
and suppose $q \in I$. 
As $q \notin (x,y)^2$, the Hilbert function of $(R/I)_{(x,y)}$ is given by $(1,1,1)$. 
Note that the coefficients $(ac+dx)$ and $(bc+ey)$ of $x$ and $y$ in $q$ have a common divisor $c$ modulo $(x,y)$. 
Unlike the multiplicity 2 case, this does not lead to a contradiction, but merely implies that $(ax+by)x$ and $(ax+by)y$ must belong to $I$ as well, for $q(x,y) \equiv c(ax+by)(x,y)$ modulo $(x,y)^3$ and $c$ is a non-zerodivisor modulo $I$. 
Indeed, 
\begin{align*}
	I	& = (x,y)^3 + (ax+by)(x,y) + (q) \\
		& = \big( x^3,\, x^2y,\, xy^2,\, y^3,\, ax^2+bxy,\, axy+by^2,\, acx+bcy+dx^2+ey^2 \big)
\end{align*}
is an ideal of multiplicity 3 and primary to $(x,y)$. 
\end{example*}

As a partial result, we characterize in Lemma~\ref{lemma:triple-colon-linear} all $(x,y)$-primary ideals $I$ of multiplicity 3 whenever $I : (x,y)$ contains a linear form $l$. 
(Note that after choosing suitable generators for the ideal $(x,y)$, we may assume $l=x$.) 
We will require the following lemma which, similarly to Lemma~\ref{lemma:unmix}, describes yet another class of unmixed ideals. 

\begin{lemma}
\label{lemma:more-unmix}
Let $I=(x^2,xy,y^2v,cx+dyv)$ with linear forms $x,y,v$ such that $\codim{x,yv}=2$ and forms $c,d \in R$ such that $\deg(c)=\deg(d)+1$ and $(c,d) \not\subset (x,y)$. 
Then $\pd{R/I} \leq 3$. 
Further, $I$ is unmixed if and only if $\codim{x,y,c,d}>3$, in which case $R/I$ has multiplicity 3. 
\end{lemma}
(The hypothesis $\codim{x,yv}=2$ merely says that $x,y$ as well as $x,v$ are independent linear forms. 
In particular, $\codim{I}=2$.) 
\begin{proof}
We proceed as in the proof of Lemma~\ref{lemma:unmix}. 
Consider the complex 
\begin{multline*}
	0 \to R
	\xrightarrow{ \varphi_3=\left(\begin{smallmatrix} c \\ d \\ y \\ x \end{smallmatrix}\right) } R^4
	\xrightarrow{ \varphi_2=\left(\begin{smallmatrix} -y & 0 & c & 0 \\ x & -yv & dv & -c \\ 0 & x & 0 & -d \\ 0 & 0 & -x & y \end{smallmatrix}\right) } R^4 
	\xrightarrow{ \varphi_1=(\begin{smallmatrix} x^2 & xy & y^2v & cx+dyv \end{smallmatrix}) } 
	R \to 0.
\end{multline*}
We have $I_1(\varphi_1)=I$ and $I_3(\varphi_2)=I\,(x,y,c,d)$, both ideals of height two, and $I_1(\varphi_3)=(x,y,c,d)$ of height at least 3. 
So the above complex resolves $R/I$ and $\codim{P} \leq 3$ for any prime ideal $P \in \ass{R/I}$. 
If $\codim{x,y,c,d}>3$, then $I_1(\varphi_3) \not\subset P$ while $I_2(\varphi_3)=0$. 
In this case $\varphi_3$ splits locally at $P$ and $\codim{P}=\pd{R_P/I_P} = 2$.
That is, $I$ is unmixed.

Conversely, assume $c$ and $d$ have a common divisor $e \in \m$ modulo $(x,y)$. 
Write $c \equiv c'e$ and $d \equiv d'e$ modulo $(x,y)$ and note that the condition $(c,d) \not\subset (x,y)$ implies $e \notin (x,y)$, that is, $\codim{x,y,e}=3$. 
Reducing $cx+dyv$ modulo $(x^2,xy,y^2v)$, we obtain $e(c'x+d'yv) \in I$ and therefore $(x,y,e) \subseteq I : (c'x+d'yv)$. 
As $c'x+d'yv \notin I$, this means that $x,y,e$ are contained in some associated prime of $R/I$. 
As $I$ has height two, it cannot be unmixed. 

Lastly, assuming $I$ is unmixed, we determine its associated primes and compute its multiplicity using the associativity formula~\eqref{eqn:assoc}. 
Since $x^2, y^2v \in I$, any prime ideal containing $I$ must contain $x$ and either $y$ or $v$. 
In particular, we have $(x,y), (x,v) \in \ass{R/I}$.
But if $I$ is unmixed, then these are the only associated primes of $R/I$. 
If $v \in (x,y)$, then $I$ is primary to $(x,y)$ and $e(R/I) = \lambda\!\left( (R/I)_{(x,y)} \right) = 3$. 
And if $v \notin (x,y)$, then $e(R/I) = \lambda\!\left( (R/I)_{(x,y)} \right) + \lambda\!\left( (R/I)_{(x,v)} \right) = 2 + 1 = 3$. 
\end{proof}

\begin{lemma}
\label{lemma:triple-colon-linear}
If $I$ is an ideal of multiplicity 3, primary to $(x,y)$ with independent linear forms $x,y$, and $x \in I:(x,y)$, then either $I=(x,y)^2$ or $I=(x^2,xy,y^3,cx+dy^2)$ with forms $c$ and $d$ such that $\codim{x,y,c,d}>3$. 
\end{lemma}
\begin{proof}
As $\lambda\!\left( (R/I)_{(x,y)} \right) = e(R/I) = 3$, we must have $(x,y)^3 \subsetneq I$. 
And by our hypothesis $x(x,y) \subset I$. 
Thus, $(x^2,xy,y^3) \subsetneq I$. 
In addition, as $e(R/I)=3$, $I$ must contain terms of the form $c_ix+b_iy+d_iy^2$ with $c_i \notin (x,y)$. 
Assuming $I \neq (x,y)^2$, we may choose $b_i=0$. 
Indeed, multiplying $c_ix+b_iy+d_iy^2$ with $y$ yields $b_i y^2 \in I$. 
As $I \neq P^2$, we have $y^2 \notin I$ and so $b_i \in (x,y)$. 
After reducing $b_iy$ modulo $xy$ and relabeling $d_i$, we can rewrite $c_ix+b_iy+d_iy^2$ as $c_ix+d_iy^2$. 

It now follows that $I$ contains only one such term as a minimal generator, that is, $I=(x^2,xy,y^3,cx+dy^2)$. 
Indeed, by the same argument as in the proof of Proposition~\ref{prop:class}, mutatis mutandis, if $cx+dy^2$ is a minimal generator of $I$, then $\codim{x,y,c,d}>3$. 
Thus, by Lemma~\ref{lemma:more-unmix}, the ideal $(x^2,xy,y^3,cx+dy^2)$ is unmixed and by Lemma~\ref{lemma:unmix-eq} it equals $I$. 
\end{proof}

\section{Projective dimension of three cubics}
\label{section:bounds}

Let $f,g,h \in R$ be three cubic forms and let $I=\fgh\unm$ denote the unmixed part of the ideal $\fgh$. 
In this section we prove that if $I$ contains a quadric, then $\pdfgh \leq 4$. 

By Remark~\ref{rmk:codim-2} we may assume $\codim{f,g,h}=2$. 
For a height two ideal, the assumption that $I$ contains a quadric implies that its multiplicity is bounded above by 6. 
For the cases with multiplicity 1 or 2, we simply draw on results from Section~\ref{section:low-mult} to prove our claim directly. 
For the cases with multiplicity $\geq 4$ we apply linkage; 
we choose a quadric $q \in I$ and a cubic $p \in I$ such that $q$ and $p$ form a regular sequence and we consider the link $I'=(q,p):I$ which has multiplicity $6-e(R/I) \leq 2$. 
Due to its low multiplicity, we are able to bound the projective dimension of the link $R/I'$ and apply Theorem~\ref{thm:extend} to prove our claim. 

This leaves us with the borderline case of multiplicity 3 which will require most of our attention. 
To prepare for that case, we first go through a relatively lengthy and at times technical analysis of all height two unmixed ideals $I$ of multiplicity 3, and consider these as the unmixed part of $\fgh$. 
For the most part, we either show that $R/I$ is \cm{} or we bound the projective dimension of an ideal linked to $I$. 
In the process, we establish the following fact.

\begin{prop}
Let $R$ be a polynomial ring over a field and let $f,g,h \in R$ be three cubic forms. 
If $R/\fgh$ has multiplicity 3, then $\pdfgh \leq 16$. 
\end{prop}

\medskip

So suppose $\codim{f,g,h}=2$, $e(R/\fgh)=3$, and let $I=\fgh\unm$. 
By the associativity formula~\eqref{eqn:assoc} there are five cases to consider: 
\begin{enumerate}
	\item $I$ is a prime ideal of multiplicity 3. 
	\item $I$ is primary to a prime ideal $P$ of multiplicity 1 and $(R/I)_P$ has length 3. 
	\item $I$ is the intersection of two prime ideals with multiplicities 1 and 2, respectively. 
	\item $I$ is the intersection of a prime ideal of multiplicity 1 and a primary ideal of multiplicity 2. The latter ideal is of the form as described in part~(iv) or (iv$^\circ$) of Proposition~\ref{prop:class}. 
	\item $I$ is the intersection of three prime ideals, each of multiplicity 1. 
\end{enumerate}

\medskip

\noindent
\underline{Case 1}~ 
$I$ is a prime ideal of multiplicity 3. 
By Proposition~\ref{prop:degen} we may assume $I$ is non-degenerate. 
Thus, $I$ is a homogeneous prime ideal of \emph{minimal multiplicity}, that is, $e(R/I) = \codim{I} + 1$. 
It is well known that $R/I$ is \cm{} (cf.~\cite{eh}) and we have $\pdfgh \leq 3$ by Theorem~\ref{thm:fan}. 

\medskip

\noindent
\underline{Case 2}~ 
$I$ is primary to $P=(x,y)$ with independent linear forms $x,y$ and the Hilbert function of $(R/I)_P$ is either $(1,2)$ or $(1,1,1)$. 
That is, locally at $P$, the ideal $I$ is either of the form $(x,y)_P^2$ or of the form $(x,y)_P^3 + (cx+dy)^{}_P$ with elements $c,d$ such that $(c,d) \not\subset (x,y)$. 

In the former case we must have $I \subseteq P^2$, as otherwise $(R/I)_P$ would have Hilbert function $(1,1,1)$. 
But $\ass{R/I}=\{P\}$ and $I^{}_P=P_P^2$. 
So $I=P^2$ globally and $\pdfgh \leq 3$ by Theorem~\ref{thm:fan}. 

Now suppose $(R/I)_P$ has Hilbert function $(1,1,1)$. 
Note that $I:P$ is also primary to $P$ and has multiplicity 2. 
By parts~(iv) or (iv$^\circ$) of Proposition~\ref{prop:class}, either $I:P = (x,y)^2 + (ax+by)$ where $x,y,a,b$ form a regular sequence, or $I:P = (x,y^2)$. 

If $I:P = (x,y^2)$, then the mere fact that $I:P$ contains a linear form allows us to give an explicit description of $I$ in terms of its generators. 
By Lemma~\ref{lemma:triple-colon-linear}, $I=(x^2,xy,y^3,cx+dy^2)$ with elements $c,d$ such that $\codim{x,y,c,d}>3$. 
Having determined $I$ explicitly, we are able to compute its link $(x^2,y^3):I = (x^2,xy,y^3,cx-dy^2)$. 
Now, by Lemma~\ref{lemma:more-unmix}, $\pd{R/(x^2,xy,y^3,cx-dy^2)} \leq 3$ and so $\pdfgh \leq 4$ by Theorem~\ref{thm:extend}. 

For the remaining case where $I:P = P^2 + (ax+by)$, we first point out that $\deg(ax+by) \leq 3$. 
Indeed, we have 
\[
	\fgh ~ \subseteq ~ I ~ \subset ~ I:P ~ = ~ P^2 + (ax+by).
\]
So if $\deg(ax+by) \geq 4$, then $\fgh \subset P^2$ and by Lemma~\ref{lemma:unmixed-hull}, $I \subseteq P^2$ --- a contradiction, as $(R/I)_P$ has Hilbert function $(1,1,1)$. 

If $\deg(ax+by) = 3$, then there are linear forms $l_{ij}$ and field coefficients $\alpha, \beta, \gamma$ such that 
\[
	\begin{pmatrix}f \\ g \\ h\end{pmatrix}
	=
	\begin{pmatrix}
		l_{11} & l_{12} & l_{13} & \alpha \\
		l_{21} & l_{22} & l_{23} & \beta \\
		l_{31} & l_{32} & l_{33} & \gamma \\
	\end{pmatrix}
	\begin{pmatrix}x^2 \\ xy \\ y^2 \\ ax+by\end{pmatrix}.
\]
As $\fgh \not\subset (x,y)^2$, one of the coefficients $\alpha, \beta, \gamma$ is non-zero; 
say $\alpha \neq 0$.
Setting $a' = l_{11}x+l_{12}y+\alpha a$ and $b' = l_{13}y+\alpha b$, we have $f=a'x+b'y$.
As $\alpha \neq 0$ and the elements $x,y,a,b$ form a regular sequence, so do the elements $x,y,a',b'$.
By Lemma~\ref{lemma:unmix} the ideal $P^3+(f)$ is unmixed of multiplicity 3 and by Lemma~\ref{lemma:unmix-eq} it is equal to $I$. 
This allows us now to compute the link $(x^3,y^3) : I = \left(x^3, x^2y^2, y^3, (a'x-b'y)xy, a'^2x^2-a'b'xy+b'^2y^2\right)$, whose quotient has projective dimension 3. 
Thus, $\pdfgh \leq 4$ by Theorem~\ref{thm:extend}. 

It remains the case $\deg(ax+by) = 2$. 
The cubics $f,g,h$ can be expressed in terms of the quadrics $x^2, xy, y^2, ax+by$ using no more than 12 linear forms $l_{ij}$. 
(Note that $a,b$ are linear as well.) 
Without having determined the unmixed part $I$ in this case, we use the fact $f,g,h \in k[l_{ij},a,b,x,y]$ to infer that $\pdfgh \leq 16$. 

\medskip

\noindent
\underline{Case 3}~ 
$I$ is the intersection of two height two prime ideals with multiplicities 1 and 2, respectively.
Write $I = (u,v) \cap (x,q)$ with linear forms $u,v,x$ and an irreducible quadric $q$. 
After subtracting a suitable multiple of $x$ from $q$, we may assume $q$ is reduced modulo $x$ without changing the ideal $(x,q)$. 
As there is no containment among the ideals $(u,v)$ and $(x,q)$, we have $\codim{u,v,x,q} = 3$ or $4$. 

If $\codim{u,v,x,q} = 3$, then either $x$ or $q$ must belong to $(u,v)$. 
Indeed, as $q$ is reduced modulo $x$, the condition $q \notin (u,v)$ is tantamount to $q \notin (u,v,x)$. 
So if in addition $x \notin (u,v)$, then $u,v,x,q$ would form a regular sequence and $\codim{u,v,x,q} = 4$. 
Thus, $\pd{R/(u,v,x,q)} = 3$ as $(u,v,x,q)$ is generated by a regular sequence of length three --- either by $(u,v,q)$ if $x \in (u,v)$ or by $(u,v,x)$ if $q \in (u,v)$. 
Applying Lemma~\ref{lemma:depth} to the short exact sequence 
\[
	0 \longrightarrow \frac{R}{I} \longrightarrow \frac{R}{(u,v)} \oplus \frac{R}{(x,q)} \longrightarrow \frac{R}{(u,v,x,q)} \longrightarrow 0
\]
yields $\pd{R/I} = 2$ and $\pdfgh \leq 3$ by Theorem~\ref{thm:fan}. 

And if $\codim{u,v,x,q} = 4$, then $u,v,x,q$ form a regular sequence and $I=(ux,uq,vx,vq)$. 
We compute the link $(ux,vq) : I = (x,v) \cap (u,q)$, whose quotient has projective dimension 3. 
Thus, $\pdfgh \leq 4$ by Theorem~\ref{thm:extend}.

\medskip

\noindent
\underline{Case 4}~ 
By Proposition~\ref{prop:class}, $I$ admits a primary decomposition either of the form $(u,v) \cap (x,y^2)$ or of the form $(u,v) \cap (x^2,xy,y^2,ax+by)$ with independent linear forms $u,v$, independent linear forms $x,y$, and elements $a,b$ such that $x,y,a,b$ form a regular sequence. 
Since there is no containment among the associated primes $(x,y)$ and $(u,v)$, we have $\codim{x,y,u,v} = 3$ or $4$. 

The case $I=(u,v) \cap (x,y^2)$ is entirely analogous to case 3 above, with the quadric $q$ replaced by $y^2$. 
(The arguments used did not rely on the fact that $q$ was irreducible.) 
So we have $I = (u,v) \cap (x^2,xy,y^2,ax+by)$ where $x,y,a,b$ form a regular sequence. 

Removing any multiples of $x$ or $y$ from $a$ and $b$ does not change the ideal $(x^2,xy,y^2,ax+by)$, as this amounts to the reduction of the term $ax+by$ modulo $(x,y)^2$. 
Hence we may assume $a$ and $b$ are reduced modulo $(x,y)$.

Throughout the subsequent arguments we will use the following simple fact in order to analyze the intersection $(u,v) \cap (x^2,xy,y^2,ax+by)$. 

\begin{lemma}
\label{lemma:intersect}
If $K_1, K_2, L$ are ideals with $K_2 \subseteq L$, then $L \cap (K_1 + K_2) = ( L \cap K_1 ) + K_2$. 
\end{lemma}
\begin{proof}
The inclusion ``$\supseteq$'' is clear. 
As for ``$\subseteq$'', let $l = k_1+k_2$ be an element in $L \cap (K_1 + K_2)$ with $l \in L$, $k_1 \in K_1$, and $k_2 \in K_2$. 
By assumption $k_2 \in L$. So $k_1 = l - k_2 \in L \cap K_1$ and $k_1 + k_2 \in ( L \cap K_1 ) + K_2$. 
\end{proof}

And as argued previously in case 2, we have $\deg(ax+by) \leq 3$. 
For if $\deg(ax+by) \geq 4$, then the inclusion $\fgh \subseteq I = (u,v) \cap (x^2,xy,y^2,ax+by)$ would imply $\fgh \subseteq (u,v) \cap (x,y)^2$ and by Lemma~\ref{lemma:unmixed-hull}, $I \subseteq (u,v) \cap (x,y)^2$ --- a contradiction. 
In what follows we will consider the following cases individually: 
\[\begin{array}{ll}
\mathrm{Case~4.1} & \deg(ax+by)=1 \\
\mathrm{Case~4.2~(a.i)} & \deg(ax+by)=2,\; \codim{x,y,u,v} = 3,\; ax+by \in (u,v) \\
\mathrm{Case~4.2~(a.ii)} & \deg(ax+by)=2,\; \codim{x,y,u,v} = 3,\; ax+by \notin (u,v) \\
\mathrm{Case~4.2~(b.i)} & \deg(ax+by)=2,\; \codim{x,y,u,v} = 4,\; ax+by \in (u,v) \\
\mathrm{Case~4.2~(b.ii)} & \deg(ax+by)=2,\; \codim{x,y,u,v} = 4,\; ax+by \notin (u,v) \\
\mathrm{Case~4.3~(a)} & \deg(ax+by)=3,\; \codim{x,y,u,v} = 3 \\
\mathrm{Case~4.3~(b)} & \deg(ax+by)=3,\; \codim{x,y,u,v} = 4 \\
\end{array}\]

\smallskip

\noindent
\underline{Case 4.1}~ 
The case $\deg(ax+by)=1$ is equivalent to $I=(u,v) \cap (x,y^2)$ which was discussed above. 

\smallskip

\noindent
\underline{Case 4.2 (a.i)}~ 
We have $\deg(ax+by)=2$, that is, $x,y,a,b$ are independent linear forms. 
In particular, after choosing suitable generators for the ideal $(x,y)$, we may assume that $a$ is reduced modulo $b$ (or vice versa). 
More concretely, if, say, $a=a'+\beta b$, then $ax+by = a'x+b(\beta x+y)$ and we can relabel $a'$ as $a$ and $\beta x+y$ as $y$ without changing the ideal $(x,y)^2+(ax+by)$. 

As $\codim{x,y,u,v} = 3$, we may assume $u=x$. 
We also know that $y \notin (x,v)$. 
Since $ax+by \in (x,v)$, Lemma~\ref{lemma:intersect} implies 
\[
	I ~ = ~ 
	\left[(x,v) \cap (x,y)^2\right] + (ax+by) ~ = ~ 
	(x^2, xy, y^2v, ax+by). 
\]
However, $ax+by \in (x,v)$ is equivalent to $by \in (x,v)$. 
And as $y \notin (x,v)$ and $b$ is reduced modulo $x$, we have $b \in (v)$. 
Since $b$ is linear, we may assume $b=v$. 
One now sees that the generator $y^2v$ is redundant and $I$ is generated by the $2 \times 2$ minors of a $3 \times 2$ matrix: 
\[
	I ~ = ~ (x^2, xy, ax+by) ~ = ~ 
	I_2 \! \begin{pmatrix}
	x & 0 \\
	a & -y \\
	b & x
	\end{pmatrix}.
\]
So $R/I$ is \cm{} and $\pdfgh \leq 3$ by Theorem~\ref{thm:fan}. 

\smallskip

\noindent
\underline{Case 4.2 (a.ii)}~ 
As in the preceding case, we have $\deg(ax+by)=2$ and $u=x$, but now $ax+by \notin (x,v)$. 
We claim that $I = \left(x^2, xy, y^2v, (ax+by)v\right)$. 
Indeed, the ideal $\left(x^2, xy, y^2v, (ax+by)v\right)$ has height two and it is contained in $I=(x,v) \cap (x^2,xy,y^2,ax+by)$. 
So, by Lemma~\ref{lemma:unmix-eq}, it suffices to show that it is unmixed of multiplicity 3. 
This in turn will follow from Lemma~\ref{lemma:more-unmix} once we verify that $\codim{x,y,av,b}=4$. 
As $x,y,a,b$ already form a regular sequence, it suffices to show that $b \notin (x,y,v)$. 
Given that $b$ is reduced modulo $(x,y)$, this is tantamount to $b \notin (v)$ which follows from the hypothesis $ax+by \notin (x,v)$. 
So $\codim{x,y,av,b}=4$ and $I = \left(x^2, xy, y^2v, (ax+by)v\right)$. 
We now claim that 
\[
	(x^2,y^2v) : I ~ = ~ \left( x^2, xy, y^2v, (ax-by)v \right).
\]
To show this, we first observe that $\left( x^2, xy, y^2v, (ax-by)v \right) \subseteq (x^2,y^2v):I$ which amounts to checking $\big(xy, (ax-by)v\big) \big(xy, (ax+by)v\big) \subseteq (x^2,y^2v)$. 
Next we observe that $(x^2,y^2v):I$ has height two and multiplicity 3, and by Lemma~\ref{lemma:more-unmix} so does the ideal $\left( x^2, xy, y^2v, (ax-by)v \right)$. 
Our claim now follows from Lemma~\ref{lemma:unmix-eq}. 
Lemma~\ref{lemma:more-unmix} also asserts that $\pd{\sfrac{R}{(x^2,y^2v):I}} \leq 3$ and thus $\pdfgh \leq 4$ by Theorem~\ref{thm:extend}. 

\smallskip

\noindent
\underline{Case 4.2 (b.i)}~ 
As $ax+by \in (u,v)$, we have $I = \left[ (u,v) \cap (x,y)^2 \right] + (ax+by)$ by Lemma~\ref{lemma:intersect}. 
In addition, $\codim{x,y,u,v}=4$ implies $a \in (u,v,y)$ and $b \in (u,v,x)$. 
Since $a,b$ are reduced modulo $(x,y)$, we have $a,b \in (u,v)$. 
But $a,b$ are linear and linearly independent, so $(u,v)=(a,b)$ and we obtain 
\[
	I ~ = ~ (ax^2, bx^2, ay^2, by^2, ax+by). 
\]
The link $(ax^2, by^2):I$ is the ideal $\left(ax^2, by^2, x^2y^2, abxy, (ax-by)ab\right)$, whose quotient has projective dimension 3. 
Hence $\pdfgh \leq 4$ by Theorem~\ref{thm:extend}. 

\smallskip

\noindent
\underline{Case 4.2 (b.ii)}~ 
We have $ax+by \notin (u,v)$ which implies $(a,b) \not\subset (u,v)$. 
Say $a \notin (u,v)$. 
Moreover, as $a$ is reduced modulo $(x,y)$, we have $a \notin (x,y,u,v)$ and $\codim{x,y,u,v,a} = 5$. 
If $b \notin (x,y,u,v,a)$, then $x,y,u,v,a,b$ are independent linear forms and 
\[
	I ~ = ~ (u,v) \, \left( (x,y)^2 + (ax+by) \right). 
\]
The link $(x^2u, y^2v) : I$ is the ideal $\left( x^2u, y^2v, x^2y^2, xyuv, (ax-by)uv \right)$, whose quotient has projective dimension 3. 
By Theorem~\ref{thm:extend}, $\pdfgh \leq 4$. 

If on the other hand $b \in (x,y,a,u,v)$, then, after reduction modulo $(x,y,a)$, we have $b = \alpha u + \beta v$. 
Without loss of generality $\alpha \neq 0$ and we may replace $u$ by $b$. 
Now $x,y,a,b,v$ are independent linear forms and as above 
\[
	I ~ = ~ (b,v) \, \left( (x,y)^2 + (ax+by) \right).
\]
The link $(x^2b, y^2v) : I$ is the ideal $\left( x^2b, y^2v, x^2y^2, xybv, (ax-by)bv \right)$, whose quotient has projective dimension 3. 
By Theorem~\ref{thm:extend}, $\pdfgh \leq 4$. 

\medskip

We are left with the cases 4.3 (a) and 4.3 (b), where $a,b$ are quadrics and $\deg(ax+by)=3$. 
We first reduce to the situation where $ax+by \in (u,v)$. 
Recall that we cannot have $\fgh \subset (x,y)^2$, for otherwise $I \subset (u,v) \cap (x,y)^2$ by Lemma~\ref{lemma:unmixed-hull}. 
So suppose, without loss of generality, that $f$ has a non-zero contribution from the cubic term $ax+by$, that is, 
\[
	f ~ = ~ l_1 \, x^2 + l_2 \, xy + l_3 \, y^2 + \alpha (ax+by)
\]
with linear forms $l_1, l_2, l_3$ and a scalar $0 \neq \alpha \in k$. 
Setting $a' = l_1 \, x+l_2 \, y+\alpha a$ and $b' = l_3 \, y+\alpha b$, we have $f=a'x+b'y$ where $x,y,a',b'$ form a regular sequence. 
By Lemma~\ref{lemma:unmix} the ideal $(x,y)^2+(f)$ is unmixed of multiplicity 2 and by Lemma~\ref{lemma:unmix-eq} it is equal to $(x,y)^2+(ax+by)$. 
So we may replace $ax+by$ by $f=a'x+b'y$ without changing the ideal $I$. 
Note that $f \in I \subset (u,v)$. 
To ease notation, we relabel $a'$ as $a$ and $b'$ as $b$ and arrive at $I = \left[(u,v) \cap (x,y)^2\right] + (ax+by)$ with $ax+by \in (u,v)$. 

\smallskip

\noindent
\underline{Case 4.3 (a)}~ 
With $\codim{x,y,u,v} = 3$ we may assume $u=x$. 
By Lemma~\ref{lemma:intersect}, $I=(x^2,xy,y^2v,ax+by)$. 
Note that $y \notin (x,v)$. 
So $ax+by \in (x,v)$ implies $b \in (x,v)$, say $b=b_1x+b_2v$ with linear forms $b_1,b_2$. 
As $x,y,a,b$ form a regular sequence, so do $x,y,a,b_2$. 
We have $ax+by = (a+b_1y)x+b_2vy$. 
Relabeling $a+b_1y$ as $a$ and $b_2$ as $b'$, we arrive at $I=(x^2,xy,y^2v,ax+b'yv)$ with $\codim{x,y,a,b'}=4$. 
We now claim that 
\[
	(x^2, y^2v) : I ~ = ~ (x^2, xy, y^2v, ax-b'vy). 
\]
This is shown by the exact same arguments as carried out in case 4.2~(a.ii). 
By Lemma~\ref{lemma:more-unmix}, $\pd{\sfrac{R}{(x^2, y^2v):I}} \leq 3$ and so $\pdfgh \leq 4$ by Theorem~\ref{thm:extend}. 

\smallskip

\noindent
\underline{Case 4.3 (b)}~ 
With $\codim{x,y,u,v} = 4$ we have 
\begin{equation}
\label{not-any-quadrics}
	I ~ = ~ (x^2u,x^2v,xyu,xyv,y^2u,y^2v,ax+by). 
\end{equation}
(Note that in this case $I$ is generated by cubics and does not contain any quadrics.) 
As $x,y,u,v$ form a regular sequence, $ax+by \in (u,v)$ implies $a \in (y,u,v)$ and $b \in (x,u,v)$. 
A priori we would need six linear coefficients to express the quadrics $a$ and $b$ in terms of the linear forms $x,y,u,v$. 
However, after combining the coefficients of $xy$ in $ax+by$, this number can be reduced to 5. 
So, the generators of $I$ can be expressed entirely in terms of 9 linear forms. 
As $\fgh \subseteq I$, we have $\pdfgh \leq 9$. 

\medskip

\noindent
\underline{Case 5}~ 
$I$ is the intersection of three prime ideals, each of height two and multiplicity one. 
Write $I=(x,y)\cap(u,v)\cap(s,t)$ with linear forms $x,y,u,v,s,t$ and note that $3 \leq \codim{x,y,u,v,s,t} \leq 6$. We consider each case by giving an explicit description of the ideal $I$, up to a linear change of coordinates. 

\smallskip

In the height three case, $I = (x,y) \cap (y,u) \cap (u,x) = (xy,yu,ux)$. 
$R/I$ is \cm{} and $\pdfgh \leq 3$ by Theorem~\ref{thm:fan}.

\smallskip

In the height four case, either $I = (x,y) \cap (y,u) \cap (u,v) = (xu,yu,yv)$ or $I = (x,y) \cap (x,u) \cap (x,v) = (x,yuv)$. 
In both cases $R/I$ is \cm{} and $\pdfgh \leq 3$. 

\smallskip

In the height five case, $I = (x,y) \cap (u,v) \cap (v,s) = (xv,yv,xus,yus)$. 
The link $(xv,yus) : I$ is the ideal $(xy,xv,yus,vus)$, whose quotient has projective dimension 3. 
Hence $\pdfgh \leq 4$ by Theorem~\ref{thm:extend}. 

\smallskip

In the height six case, $I=(xus,xut,xvs,xvt,yus,yut,yvs,yvt)$. 
The link $(xus,yvt) : I$ is the ideal $(xus,yvt,xyuv,xyst,uvst)$, whose quotient has projective dimension 3. 
Again, $\pdfgh \leq 4$. 

\bigskip

We are now ready to prove the main result. 

\begin{thm}
\label{thm:unmix-quad}
Let $R$ be a polynomial ring over a field and let $J \subset R$ be an ideal generated by three cubic forms. 
If the unmixed part of $J$ contains a quadric, then $\pd{R/J} \leq 4$.
\end{thm}
\begin{proof}
By Remark~\ref{rmk:codim-2} we may assume $\codim{J}=2$. 
Furthermore, for a height two, three-generated ideal $J$ we showed in Section~\ref{section:low-mult}, without any assumptions on the unmixed part or the degrees of the generators, that $\pd{R/J} \leq 4$ if $e(R/J)\leq2$. 
Thus, we may further assume $e(R/J) \geq 3$.

Now we utilize the assumption that the unmixed part of $J$ contains a quadric. 
Denote by $I$ the unmixed part of $J$ and let $q \in I$ be a quadric. 
As $J$ is generated by cubics and $\codim{J}=2$, we can invoke prime avoidance and choose a cubic form $p \in J \subseteq I$ such that $q$ and $p$ form a regular sequence. 
This imposes an upper bound on the multiplicity of $R/I$, namely, $e(R/I) \leq e(R/(q,p)) = 6$. 
By Lemma~\ref{lemma:unmix-eq} equality holds if and only if $I=(q,p)$, in which case $\pd{R/J} \leq 3$ by Theorem~\ref{thm:fan}. 
So there remain the cases $3 \leq e(R/I) \leq 5$ to consider. 

Consider the link $(q,p):I$ which has multiplicity $6-e(R/I)$. 
If $e(R/I)=5$, then the link $(q,p):I$ is a height two unmixed ideal of multiplicity 1, so it is generated by two linear forms and $\pd{R/J} \leq 3$ by Theorem~\ref{thm:extend}. 
If $e(R/I)=4$, then $(q,p):I$ is a height two unmixed ideal of multiplicity 2. 
Such ideals were classified in Proposition~\ref{prop:class} whereby $\pd{\sfrac{R}{(q,p):I}} \leq 3$. 
Thus, $\pd{R/J} \leq 4$ by Theorem~\ref{thm:extend}.

It remains the case $e(R/I)=3$. 
Note that now the link $(q,p):I$ has multiplicity 3 as well. 
Now we resort to the arguments carried out earlier in this section: 
If $I$ is an ideal as described in cases 1, 3, or 5, then we have already shown (without the assumption that $I$ contains a quadric) that $\pd{R/J} \leq 4$. 
As for ideals in case 4, the only instance where we obtained a bound for $\pd{R/J}$ greater than 4 was that of the unmixed part generated entirely by cubics --- cf.~\eqref{not-any-quadrics}. 
So we are left with the remaining case 2, that is, $I$ is primary to $(x,y)$ with independent linear forms $x,y$. 
Note that $(x,y)^3 \subsetneq I$.

First we assume that $q \notin (x,y)^2$. 
Say $q=cx+dy$ with $(c,d) \not\subset (x,y)$, that is, $\codim{x,y,c,d} \geq 3$. 
If $\codim{x,y,c,d}=4$, then the ideal $(x,y)^3+(cx+dy)$ is unmixed by Lemma~\ref{lemma:unmix} and $I = (x,y)^3+(cx+dy)$ by Lemma~\ref{lemma:unmix-eq}. 
In this case the link $(x^3,y^3) : I = \left(x^3, x^2y^2, y^3, (cx-dy)xy, x^2c^2-xycd+y^2d^2\right)$ has projective dimension 3 and $\pd{R/J} \leq 4$ by Theorem~\ref{thm:extend}. 

If on the other hand $\codim{x,y,c,d}=3$, then we may choose $q$ to be of the form $cx+\alpha xy+y^2$ with some field coefficient $\alpha$. 
Indeed, suppose $\codim{x,y,c,d}=3$, say, $c \notin (x,y)$ and $d=\alpha x+\beta y+\gamma c$ with field coefficients $\alpha, \beta, \gamma$. 
We point out that $\beta \neq \alpha \gamma$, for if $\beta=\alpha\gamma$, then $cx+dy=(c+\alpha y)(x+\gamma y)$ and $c+\alpha y$ would be a zerodivisor on $I$, contradicting $c \notin (x,y)$. 
We have 
\begin{align*}
	q	& = cx+(\alpha x+\beta y+\gamma c)y \\
		& = c(x+\gamma y)+(\alpha x+\beta y)y \\
		& = c(x+\gamma y)+(\alpha (x+\gamma y)+(\beta-\alpha\gamma) y)y.
\end{align*}
Relabeling $x+\gamma y$ as $x$, we can write $q=cx+\alpha xy+(\beta-\alpha\gamma)y^2$. 
As $\beta-\alpha\gamma \neq 0$, we can further rescale $y$ and write $q=cx+\alpha xy+y^2$. 
(More precisely, we relabel $\sqrt{\beta-\alpha\gamma} \, y$ as $y$ and $\alpha/\sqrt{\beta-\alpha\gamma}$ as $\alpha$.) 
If $\alpha \neq 0$, then we may also rescale $x$ and write $q$ simply as $cx+xy+y^2$. 
We shall see in retrospect that the value of $\alpha$ has no bearing on the form of the unmixed part $I$, as in either case the element $xy$ is in $I$.

We are now able to compute the link $(q,y^3):I$. 
First we compute the colon $(q,y^3):\left( (q)+(x,y)^3 \right)$ by hand, using the fact that $x,y,c$ are independent linear forms: 
\begin{equation}
\label{eqn:alpha-link}
(q,y^3):\left( (q)+(x,y)^3 \right) = \begin{cases}
	~ ( c^2-y^2, cy+y^2, cx+xy+y^2 ) & \mbox{if~~} \alpha = 1, \\
	~ ( c^2, cy, cx+y^2 ) & \mbox{if~~} \alpha = 0.
\end{cases}
\end{equation}
In either case the ideal $(q,y^3):\left( (q)+(x,y)^3 \right)$ has multiplicity 3 and it contains the ideal $(q,y^3):I$, since $(q)+(x,y)^3 \subseteq I$.
By Lemma~\ref{lemma:unmix-eq} we have equality: 
\[
	(q,y^3):I ~ = ~ (q,y^3):\left( (q)+(x,y)^3 \right).
\]
It is easily seen that both ideals in \eqref{eqn:alpha-link} have \cm{} quotients, as they are generated by the $2 \times 2$ minors of a $3 \times 2$ matrix. 
Thus, $\pd{R/I} = \pd{\sfrac{R}{(q,y^3):I}} = 2$ and $\pd{R/J} \leq 3$.

(Note: Linking the ideal $(q,y^3):I$ back to $I=(q,y^3):\left( (q,y^3):I \right)$ reveals that $I = (x^2,xy,cx+y^2)$ regardless of the value of $\alpha$.) 

It remains the case where $q \in (x,y)^2$. 
By Remark~\ref{rmk:alg-closed} we may assume without loss of generality that the ground field $k$ is algebraically closed. 
Thus we can factor the quadric $q$ as $l\,l'$ where $l,l' \in (x,y) \setminus I$ are two (not necessarily independent) linear forms.
Now consider the following chain of $(x,y)$-primary ideals: 
\[
	I ~~ \subsetneq ~~ I:(x,y) ~~ \subseteq ~~ I:(l) ~~ \subseteq ~~ (x,y).
\]
We have $e(R/I)=3$ and $e(R/(x,y))=1$. 
Furthermore, $e\!\left( \sfrac{R}{I:(x,y)} \right) = 1$ if and only if $I:(x,y)=(x,y)$ and therefore $I=(x,y)^2$. 
So we may assume $e\!\left( \sfrac{R}{I:(x,y)} \right) = 2$. 
That forces $e\!\left( \sfrac{R}{I:(l)} \right)$ to equal either 2 or 1. 
In what follows we show what each of these multiplicity values entails for the structure of $I:(l)$, and subsequently for that of $I$. 
As it turns out, the results are the same in either case. 

First suppose $I:(l)$ has multiplicity 1 and therefore $I:(l)=(x,y)$. 
After a linear change of coordinates we may relabel $l$ as $x$. 
So $x^2,xy \in I$ and by Lemma~\ref{lemma:triple-colon-linear} we have $I=(x^2,xy,y^3,cx+dy^2)$ with forms $c$ and $d$ such that $\codim{x,y,c,d}>3$. 
(We are assuming $I \neq (x,y)^2$.) 
As the link $(x^2,y^3) : I = (x^2,xy,y^3,cx-dy^2)$ has projective dimension 3, we have $\pd{R/J} \leq 4$. 

Now suppose $I:(l)$ has multiplicity 2. 
(In particular, $I:(l)=I:(x,y)$.) 
By Proposition~\ref{prop:class} we have $I:(l)=(x,y)^2+(ax+by)$ with forms $a,b$ such that $\codim{x,y,a,b}>3$. 
As $q=l\,l' \in I$ and $l' \in I:(l)$, the term $ax+by$ must be linear. 
Relabeling $ax+by$ as $x$, we have $I : (x,y) = (x,y^2)$ and again $I=(x^2,xy,y^3,cx+dy^2)$ by Lemma~\ref{lemma:triple-colon-linear} and $\pd{R/J} \leq 4$ as above. 

(Note: The scenario $e\!\left( \sfrac{R}{I:(l)} \right)=1$ corresponds to the choice of $q=x^2$ and $l=l'=x$, while $e\!\left( \sfrac{R}{I:(l)} \right)=2$ corresponds to $q=xy$, $l=y$, and $l'=x$.) 
\end{proof}

We point out that Theorem~\ref{thm:unmix-quad} cannot be further improved. 
Its conclusion cannot be strengthened, as illustrated with the following example.
Let $R=k[x,y,a,b,l_1,l_2,l_3,l_4]$ and let $J = \left(l_1 x^2 + l_2 y^2,\, l_3 xy,\, l_4 (ax+by) \right)$. 
Then $J\unm=(x,y)^2+(ax+by)$ is generated by quadrics and $\pdr{R}{R/J}=4$. 
And the hypothesis of Theorem~\ref{thm:unmix-quad} cannot be weakened. In \cite[Section 3.3]{engheta-thesis} the author has constructed an ideal $J$ generated by three cubic forms with $\pd{R/J}=5$. 

\medskip

We conclude with the following theorem which was inspired by the corresponding statement of Theorem~\ref{thm:unmix-quad}. 

\begin{thm}
Let $R$ be a polynomial ring over a field and let $J \subset R$ be an ideal generated by three cubic forms with $e(R/J) \leq 5$. 
Denote by $I$ the unmixed part of $J$ and let $I'$ be an ideal which is linked to $I$ via cubic forms. 
If $I'$ contains a quadric, then $\pd{R/J} \leq 4$. 
\end{thm}
\begin{proof}
As argued in the proof of Theorem~\ref{thm:unmix-quad}, we may assume $\codim{J}=2$ and $e(R/J) \geq 3$. 
By our hypothesis, $I'=(p_1, p_2):I$ with cubic forms $p_1, p_2 \in I$ which form a regular sequence. 
In particular, $e(R/I')=9-e(R/I)$. 

Let $q \in I'$ be a quadric. 
As $I'$ contains cubics which generate an ideal of height two (such as $p_1$ and $p_2$), we can invoke prime avoidance and choose a cubic $p \in I'$ such that $q$ and $p$ form a regular sequence. 

If $e(R/I)=3$, then $e(R/I')=6$ and $I'=(q,p)$ by Lemma~\ref{lemma:unmix-eq}. 
In this case $R/I'$ and $R/I$ are \cm{} and $\pd{R/J} \leq 3$ by Theorem~\ref{thm:fan} or by Theorem~\ref{thm:extend}. 
If $e(R/I)=4$, then $e(R/I')=5$ and we consider a further link $K=(q,p):I'$. 
As $e(R/K)=6-5=1$, $K$ is generated by two independent linear forms. 
Thus $R/K$, $R/I'$, and $R/I$ are \cm{} and again $\pd{R/J} \leq 3$. 

It remains the case $e(R/I)=5$ with $e(R/I')=4$ and $e(R/K)=6-4=2$. 
It follows from Proposition~\ref{prop:class} that $K$ contains a second quadric $q'$ such that $q$ and $q'$ form a regular sequence. 
So we consider yet another link $K'=(q,q'):K$ with $e(R/K')=4-2=2$. 
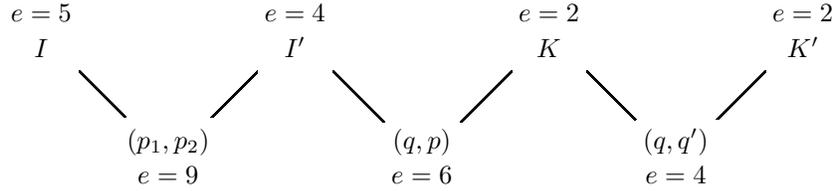
\begin{diagram}[size=2.4em]
\begin{matrix}e=5 \\ I\end{matrix} & & & & \begin{matrix}e=4 \\ I'\end{matrix} & & & & \begin{matrix}e=2 \\ K\end{matrix} & & & & \begin{matrix}e=2 \\ K'\end{matrix} \\
 & \rdLine & & \ruLine & & \rdLine & & \ruLine & & \rdLine & & \ruLine & \\
 & & \begin{matrix}(p_1,p_2) \\ e=9\end{matrix} & & & & \begin{matrix}(q,p) \\ e=6\end{matrix} & & & & \begin{matrix}(q,q') \\ e=4\end{matrix} & &
\end{diagram}
By Theorem~\ref{thm:extend}, $\pd{R/J} \leq \pd{R/I'}+1$. 
As $I'$ and $K'$ are both linked to $K$, we have $\pd{R/I'}=\pd{R/K'}$. 
And $\pd{R/K'} \leq 3$ by Proposition~\ref{prop:class}. 
Thus, $\pd{R/J} \leq 4$ as claimed. 
\end{proof}

\section*{Acknowledgments}
I thank Craig Huneke for numerous discussions, guidance, and support. 
Most computations during the preparation of this paper were performed using the computer algebra system Macaulay~2~\cite{m2}. 

\end{document}